\documentclass[11pt,reqno]{amsart}

\usepackage{amsthm,amsmath,amssymb,url,color}
\newtheorem{thm}{Theorem}[section]
\newtheorem{lem}[thm]{Lemma}
\newtheorem{cor}[thm]{Corollary}
\theoremstyle{remark}
\newtheorem*{rmk}{Remark}

\newtheorem*{eg*}{Example}

\begin{document}
\title[$q$-Faulhaber and $q$-Sali\'e coefficients]{Combinatorial
  Interpretations of the $q$-Faulhaber and $q$-Sali\'e Coefficients}

\author{Victor J. W. Guo}
\address{Institut Camille Jordan,
         Universit\'e Claude Bernard (Lyon I)\\
         F-69622, Villeurbanne Cedex, France}
\email{jwguo@eyou.com}
\author{Martin Rubey}
\address{Institut f\"ur Statistik und Decision Support\\
         Universit\"at Wien\\
         A-1010 Wien, Austria}
\email{martin.rubey@univie.ac.at}
\urladdr{http://www.mat.univie.ac.at/\textasciitilde rubey}
\author{Jiang Zeng}
\address{Institut Camille Jordan,
        Universit\'e Claude Bernard (Lyon I)\\
        F-69622, Villeurbanne Cedex, France}
\email{zeng@math.univ-lyon1.fr}

\dedicatory{Dedicated to  Xavier Viennot on the occasion of his sixtieth birthday}

\begin{abstract}
  Recently, Guo and Zeng discovered two families of polynomials featuring in a
  $q$-analogue of Faulhaber's formula for the sums of powers and a $q$-analogue
  of Gessel-Viennot's formula involving Sali\'e's coefficients for the
  alternating sums of powers. In this paper, we show that these are polynomials
  with symmetric, nonnegative integral coefficients by refining
  Gessel-Viennot's combinatorial interpretations.
\end{abstract}
\maketitle

\section{Introduction}
In the early seventeenth century, Johann Faulhaber~\cite{Faulhaber} (see also \cite{Knuth1993}) considered
the sums of powers $S_{m,n}=\sum_{k=1}^n k^m$ and provided formulas for the
coefficients $f_{m,k}$ ($0\leq m\leq 8$) in
\begin{equation}\label{eq:faul}
  S_{2m+1,n}=\frac{1}{2}\sum_{k=1}^m f_{m,k}\left(n(n+1)\right)^{k+1},
\end{equation}

In 1989, Ira Gessel and Xavier Viennot \cite{GesselViennot1989} studied the alternating sum
$T_{2m,n}=\sum_{k=1}^{n}(-1)^{n-k}k^{2m}$ and showed that there exist integers $s_{m,k}$
such that
\begin{equation}\label{eq:gv}
  T_{2m,n}=\frac{1}{2}\sum_{k=1}^{m}s_{m,k}(n(n+1))^k.
\end{equation}
In particular, they proved that the Faulhaber coefficients $f_{m,k}$ and the Sali\'e
coefficients $s_{m,k}$ count certain families of non-intersecting lattice paths.

Recently, two of the authors~\cite{GuoZeng2005}, continuing work of Michael
Schlosser~\cite{Schlosser}, Sven Ole Warnaar~\cite{Warnaar} and
Kristina Garrett and Kristen Hummel~\cite{GarrettHummel}, have found $q$-analogues of \eqref{eq:faul} and \eqref{eq:gv}.
More precisely, setting $[k]=\frac{1-q^k}{1-q}$, $[k]!=\prod_{i=1}^k[k]$,
and
  \begin{align}\label{eq:qpower}
    S_{m,n}(q) &=\sum_{k=1}^{n}\frac{[2k]}{[2]}
                               [k]^{m-1}q^{\frac{m+1}{2}(n-k)},\\
    \label{eq:altersum}
    T_{m,n}(q) &=\sum_{k=1}^{n}(-1)^{n-k}
                               [k]^{m} q^{\frac{m}{2}(n-k)},
  \end{align}
for $m,n\in\mathbb{N}$, they proved the following results:
\begin{thm}\label{thm:pqgh}
   There exist polynomials $P_{m,k}$, $Q_{m,k}$, $G_{m,k}$ and $H_{m,k}$ in
  $\mathbb{Z}[q]$ such that
  \begin{align}
    \label{eq:p}
    S_{2m+1,n}(q) &=\sum_{k=0}^m (-q^n)^{m-k} \frac{[k]!}{[m+1]!}
                                 P_{m,m-k}(q)
                                 \frac{([n][n+1])^{k+1}}{[2]},\\
    \label{eq:qmn}
    S_{2m,n}(q)   &=(1-q^{n+\frac{1}{2}})
                    \sum_{k=0}^m (- q^n)^{m-k}
                                \frac{(1-q^\frac{1}{2})^{m-k}
                                      Q_{m,m-k}(q^\frac{1}{2}) }
                                     {\prod_{i=0}^{m-k}(1-q^{m-i+\frac{1}{2}})}
                                 \frac{([n][n+1])^k}{[2]},\\
    \label{eq:t2mnq}
    T_{2m,n}(q)   &=\sum_{k=1}^{m} (-q^n)^{m-k}
                                   \frac{G_{m,m-k}(q)}
                                        {\prod_{i=0}^{m-k}(1+q^{m-i})}
                              ([n][n+1])^k,\\
    \intertext{and}
    \label{eq:t2m+1}
    T_{2m-1,n}(q) &=(-1)^{m+n}H_{m,m-1}(q^{\frac12})\frac{q^{(m-\frac{1}{2})n}}
                   {(1+q^{\frac12})^{m}\prod_{i=0}^{m-1}(1+q^{m-i-\frac12})} \nonumber\\
                  &\quad+\frac{1-q^{n+\frac12}}{1-q^{\frac12}}
                    \sum_{k=1}^{m} (-q^n)^{m-k}
                                   \frac{H_{m,m-k}(q^{\frac12})
                                        ([n][n+1])^{k-1}}
                     {(1+q^{\frac12})^{m-k+1}\prod_{i=0}^{m-k}(1+q^{m-i-\frac12})}.
  \end{align}
\end{thm}

Comparing with \eqref{eq:qpower} and \eqref{eq:altersum}, we have

\begin{align*}
  f_{m,k}&=(-1)^{m-k}\frac{k!}{(m+1)!}P_{m,m-k}(1)
\intertext{and}
  s_{m,k}&=(-1)^{m-k}2^{k-m}G_{m,m-k}(1),
\end{align*}
 but the numbers corresponding to $Q_{m,k}(1)$ and $H_{m,k}(1)$ do not seem to be studied
 in the literature.
 The first values of $P_{m,k}$, $Q_{m,k}$, $G_{m,k}$ and $H_{m,k}$ are given in
 Tables~\ref{t:p}--\ref{t:h}, respectively.

\begin{table}[h]
\caption{Values of $P_{m,k}(q)$ for $0\leq k<m\leq 5$.\label{t:p}}
{\footnotesize
\begin{center}
\begin{tabular}{|l|r|r|r|r|r|r|}
\hline
$k\setminus m$ &0 & 1 & 2 & 3 & 4 & 5\\\hline
0              & 1  & 1 & 1 & 1 & 1 & 1\\\hline
1              &&   & 1 & $2(q+1)$ & $3q^2+4q+3$ & $2(q+1)(2q^2+q+2)$\\\hline
2              &&   &   & $2(q+1)$ & $(q+1)(5q^2+8q+5)$ & $(q+1)(9q^4+19q^3+29q^2+19q+9)$\\\hline
3              &&   &   &          & $(q+1)(5q^2+8q+5)$ & $2(q+1)^2(q^2+q+1)(7q^2+11q+7)$\\\hline
4              &&   &   &          &                    & $2(q+1)^2(q^2+q+1)(7q^2+11q+7)$\\\hline
\end{tabular}
\end{center}
}
\end{table}

\begin{table}[h]
\caption{Values of $Q_{m,k}(q)$ for $0\leq k<m\leq 4$.\label{t:q}}
{\footnotesize
\begin{center}
\begin{tabular}{|l|r|r|r|r|}
\hline
$k\setminus m$ & 1 & 2 & 3 & 4 \\\hline
0              & 1 & 1 & 1 & 1 \\\hline
1              &   & 1     &  $2q^2+q+2$  & $3q^4+2q^3+4q^2+2q+3$
\\\hline
2              &   &   & $2q^2+q+2$ & $(q^2+q+1)(5q^4+q^3+9q^2+q+5)$
\\\hline
3              &   &   &          & $(q^2+q+1)(5q^4+q^3+9q^2+q+5)$ \\\hline
\end{tabular}
\end{center}
}
\end{table}

\begin{table}[h]
\caption{Values of $G_{m,k}(q)$ for $0\leq k<m\leq 5$.\label{t:g}}
{\footnotesize
\begin{center}
\begin{tabular}{|l|r|r|r|r|r|}
\hline
$k\setminus m$ & 1 & 2 & 3 & 4 & 5\\\hline
0              & 1 & 1 & 1 & 1 & 1\\\hline
1              &   & 2 & $3(q+1)$&$4(q^2+q+1)$ & $5(q+1)(q^2+1)$\\\hline
2              &   &   & $6(q+1)$&$2(q+1)(5q^2+7q+5)$&$5(q+1)(3q^4+4q^3+8q^2+4q+3)$\\\hline
3              &   &   &         &$4(q+1)(5q^2+7q+5)$&$5(q+1)^2(7q^4+14q^3+20q^2+14q+7)$\\\hline
4              &   &   &         &                   &$10(q+1)^2(7q^4+14q^3+20q^2+14q+7)$\\\hline
\end{tabular}
\end{center}
}
\end{table}

\begin{table}[h]
\caption{Values of $H_{m,k}(q)$ for $0\leq k< m\leq 4$.\label{t:h}}
{\footnotesize
\begin{center}
\begin{tabular}{|l|r|r|r|r|}
\hline
$k\setminus m$&1 & 2 &     3          & 4   \\\hline
0             &1 & 1 &     1          & 1  \\\hline
1             &  & 2 & $3q^2+2q+3$    & $4q^4+3q^3+4q^2+3q+4$ \\\hline
2             &  &   & $2(3q^2+2q+3)$ & $10q^6+15q^5+30q^4+26q^3+30q^2+15q+10$ \\\hline
3             &  &   &                & $2(10q^6+15q^5+30q^4+26q^3+30q^2+15q+10)$  \\\hline
\end{tabular}
\end{center}
}
\end{table}
Recall that a polynomial $f(x)=a_0+a_1x+\cdots+a_n x^n$ of degree $n$ has
\emph{symmetric coefficients} if $a_i=a_{n-i}$ for $0\leq i\leq n$.  The tables
above suggest that the coefficients of the polynomials $P_{m,k}$, $Q_{m,k}$,
$G_{m,k}$ and $H_{m,k}$ are nonnegative and symmetric. The aim of this paper is
to prove this fact by showing that the coefficients count certain families of
non-intersecting lattice paths.

\section{Inverses of matrices}
Recall that the $n$-\emph{th complete homogeneous functions in $r$
variables} $x_1,x_2,\ldots, x_r$ has the following generating
function:
$$
\sum_{n\geq 0}h_n(x_1,\ldots,
x_r)t^n=\frac{1}{(1-x_1t)(1-x_2t)\ldots (1-x_rt)}.
$$
For $r,s\geq 0$, let $h_{n}(\{1\}^r,\{q\}^s)$ denote the $n$-th
complete homogeneous functions in $r+s$ variables, of which $r$
are specialized to $1$ and the others to $q$, i.e.,
\begin{equation}
  \label{eq:hdef}
\sum_{n\geq0}h_{n}(\{1\}^r,\{q\}^s)z^n=\frac{1}{(1-z)^r(1-qz)^s}.
\end{equation}
By convention, $h_{n}(\{1\}^r,\{q\}^s)=0$ if $r<0$ or $s<0$. For
convenience, we also write $h_{n}(\{1,q\}^r)$ instead of
$h_{n}(\{1\}^r,\{q\}^r)$.

We first prove the following result.

\begin{lem}\label{lem:hsum} Let $a$ and $b$ be non-negative integers, then
  \begin{multline*}
    \sum_{m\geq0}\sum_{k\geq0} h_{m-2k}(\{1\}^{k+a},\{q\}^{k+b})
                               \left(\frac{q^l}{[l]^2}\right)^k z^m
    =\frac{[l]^2}{[2l]}
    \begin{cases}
       \frac{[l+1]}{[l]-[l+1]z}- \frac{q[l-1]}{[l]-q[l-1]z}
      &\text{for $a=1$, $b=1$,}\\
       \frac{1}{[l]-[l+1]z}+     \frac{q^l}{[l]-q[l-1]z}
      &\text{for $a=1$, $b=0$, }\\
       \frac{q^l}{[l]-[l+1]z}+   \frac{1}{[l]-q[l-1]z}
      &\text{for $a=0$, $b=1$.}
    \end{cases}
  \end{multline*}
\end{lem}
\begin{proof}
  Using the definition \eqref{eq:hdef} of the complete homogeneous functions we
  have
  \begin{equation*}
    \begin{split}
      &\sum_{m\geq0}\sum_{k\geq0} h_{m-2k}(\{1\}^{k+a},\{q\}^{k+b})
                               x^k z^m\\
      &=\sum_{k\geq0}\frac{x^k z^{2k}}{(1-z)^{k+a}(1-qz)^{k+b}}\\
      &=\frac{1}{(1-z)^{a-1}(1-qz)^{b-1}}
        \frac{1}{(1-z)(1-qz)-xz^2}.
\end{split}
  \end{equation*}
Setting $x=\frac{q^l}{[l]^2}$ a little calculation shows that
        the denominator of the second fraction factorizes:
        $$
      \frac{1}{(1-z)(1-qz)-xz^2}=
        \frac{1}{\left([l]-qz[l-1]\right)\left([l]-z[l+1]\right)}.
        $$
The result then follows from the standard partial fraction
decomposition.
\end{proof}

Let $X_n=\frac{[n][n+1]}{q^n}$. The following lemma might be interesting per se.
When $q=1$ it reduces to simple applications of
the binomial theorem.

\begin{lem}\label{lem:binomial}
For $k,m\geq 1$, set
\begin{align*}
c_{k,m}(q)&:=h_{2m-k}(\{1,q^2\}^{k-m+1})
          +qh_{2m-k-1}(\{1,q^2\}^{k-m+1}), \\
g_{k,m}(q)&:=h_{2m-k}(\{1\}^{k-m+1},\{q\}^{k-m})
          +h_{2m-k}(\{1\}^{k-m},\{q\}^{k-m+1}), \\
d_{k,m}(q)&:=g_{k,m}(q^2)+q g_{k-1,m-1}(q^2).
\end{align*}
For $m,l\geq 1$, we have
    \begin{align}\label{eq:diff1}
      X_l^{m+1}
     -X_{l-1}^{m+1}&=\sum_k h_{m-2k}(\{1,q\}^{k+1})
             [2l][l]^{2(m-k)}q^{-l(m-k+1)},\\
     \label{eq:inverseq}
      \frac{1-q^{l+\frac12}}{(1-q^{\frac12})q^{\frac l2}}
      X_l^m
     -\frac{1-q^{l-\frac12}}{(1-q^{\frac12})q^{\frac{l-1}{2}}}
      X_{l-1}^m
     &=\sum_k c_{m,m-k}(q^{\frac12})
             [2l][l]^{2(m-k-\frac12)}q^{-l(m-k+\frac12)},\\
    \label{eq:diff}
      X_l^m+X_{l-1}^m &=\sum_k g_{m,m-k}(q)[l]^{2(m-k)}q^{-l(m-k)},\\
      \label{eq:sumd}
      \frac{1-q^{l+\frac12}}{(1-q^\frac12)q^\frac{l}{2}}X_l^{m-1}
     +\frac{1-q^{l-\frac12}}{(1-q^\frac12)q^\frac{l-1}{2}}X_{l-1}^{m-1}
     &=\sum_k d_{m,m-k}(q^\frac12)[l]^{2(m-k-\frac12)}q^{-l(m-k-\frac12)}.
    \end{align}
\end{lem}
\begin{proof} The proof rests on the previous lemma.
  \begin{itemize}
  \item  Applying Lemma~\ref{lem:hsum} with $a=1$ and $b=1$ yields that the
    coefficient of $z^m$ in $\sum_k h_{m-2k}(\{1,q\}^{k+1})q^{lk}[l]^{-2k}$ is
    \begin{equation*}
      \frac{[l]}{[2l]}
      \left([l+1]\left(\frac{[l+1]}{[l]}\right)^m
           -q[l-1]\left(\frac{q[l-1]}{[l]}\right)^m\right).
    \end{equation*}
    Multiplying this expression with $[2l]\frac{[l]^{2m}}{q^{l(m+1)}}$ we
    obtain \eqref{eq:diff1}.

  \item Since $c_{m,m-k}(q^\frac12)=h_{m-2k}(\{1,q\}^{k+1}) +q^\frac12
    h_{m-1-2k}(\{1,q\}^{k+1})$, Equation \eqref{eq:inverseq} follows directly from the previous
    calculation.

  \item As $g_{m,m-k}(q)=h_{m-2k}(\{1\}^{k+1},\{q\}^{k}) +
    h_{m-2k}(\{1\}^k,\{q\}^{k+1})$, applying Lemma~\ref{lem:hsum} with $a=1$,
    $b=0$  and $a=0$, $b=1$,
    \begin{align*}
      &\hskip -2mm \sum_k \left(h_{m-2k}(\{1\}^{k+1},\{q\}^{k})
        + h_{m-2k}(\{1\}^k,\{q\}^{k+1}))\right)q^{lk}[l]^{-2k}\\
     & =\frac{[l]^2}{[2l]} \left(\frac{1+q^l}{[l]-[l+1]z}
        +\frac{1+q^l}{[l]-q[l-1]z} \right)
    \end{align*}
    Multiplying the coefficient of $z^m$ of this expression with
    $[l]^{2m}q^{-lm}$ we obtain \eqref{eq:diff}.

  \item Since $d_{m,m-k}(q^\frac12)=g_{m,m-k}(q)+q^\frac12g_{m-1,m-k-1}(q)$,
    Equation \eqref{eq:sumd} follows directly from the previous calculation.
  \end{itemize}
\end{proof}

The following is the main result of this section. Note that together with
Theorems~\ref{thm:pq} and \ref{thm:gh} it also provides an alternative proof of
Theorem~\ref{thm:pqgh}.
\begin{thm}\label{thm:inverse}
The inverses of the  lower triangular matrices
$$
(h_{2m-k}(\{1,q\}^{k-m+1}))_{0\leq k,m \leq n},\quad (c_{k,m}(q))_{1\leq k,m\leq n},\quad
(g_{k,m}(q))_{1\leq k,m \leq n},\quad (d_{k,m}(q))_{1\leq k,m \leq n}
$$
are respectively the lower triangular matrices
\begin{align}
&\left((-1)^{k-m}\frac{[m]!}{[k+1]!}P_{k,k-m}(q)
\right)_{0\leq k,m\leq n}, \label{eq:invp}\\
&\left((-1)^{k-m}\frac{(1-q)^{k-m+1}Q_{k,k-m}(q)}{\prod_{i=0}^{k-m}(1-q^{2k-2i+1})}
\right)_{1\leq k,m\leq n}, \label{eq:invq}\\
&\left((-1)^{k-m}\frac{G_{k,k-m}(q)}{\prod_{i=0}^{k-m}(1+q^{k-i})}
\right)_{1\leq k,m\leq n}, \label{eq:invg}\\
&\left((-1)^{k-m}\frac{H_{k,k-m}(q)}{(1+q)^{k-m+1}\prod_{i=0}^{k-m}(1+q^{2k-2i-1})}
\right)_{1\leq k,m\leq n}. \label{eq:invh}
\end{align}
\end{thm}
\begin{proof}
\begin{itemize}
\item Summing Equation
\eqref{eq:diff1} over $l$ from $1$ to $n$ and applying Equation
\eqref{eq:qpower}, we obtain
\begin{align}\label{eq:qmn1}
X_n^{m+1}
=[2]\sum_{k=0}^{\lfloor m/2 \rfloor}
h_{m-2k}(\{1,q\}^{k+1})S_{2m-2k+1,n}(q)q^{-n(m-k+1)}.
\end{align}
Plugging \eqref{eq:p} in Equation \eqref{eq:qmn1}, the right-hand side becomes
$$
\sum_{k=0}^{\lfloor m/2 \rfloor}\sum_{l=0}^{m-k}
h_{m-2k}(\{1,q\}^{k+1})(-1)^{m-k-l}\frac{[l]!}{[m-k+1]!}P_{m-k,m-k-l}(q)
X_n^{l+1}.
$$
Comparing the coefficients of $X_n^{l+1}$ we see that
$(h_{2m-k}(\{1,q\}^{k-m+1}))_{0\leq k,m \leq n}$ and \eqref{eq:invp} are indeed inverses.

\item Summing
Equation \eqref{eq:inverseq} over $l$ from $1$ to $n$ and applying
Equation \eqref{eq:qpower}, we obtain
\begin{align}
\frac{1-q^{n+\frac12}}{(1-q^{\frac12})q^{\frac n2}}
X_n^m
&=[2]\sum_{k=0}^{\lfloor m/2 \rfloor}
c_{m,m-k}(q^{\frac12})S_{2m-2k,n}(q)q^{-n(m-k+\frac12)}.
\label{eq:qmn2}
\end{align}
Substituting \eqref{eq:qmn} into \eqref{eq:qmn2} and dividing both sides by
$\frac{1-q^{n+\frac12}}{(1-q^{\frac12})q^{\frac n2}}$, we get
$$
X_n^m
=\sum_{k=0}^{\lfloor m/2 \rfloor}\sum_{l=1}^{m-k}
   c_{m,m-k}(q^{\frac12})(-1)^{m-k-l}\frac{(1-q^\frac{1}{2})^{m-k-l}
   Q_{m-k,m-k-l}(q^\frac12)}
   {\prod_{i=0}^{m-k-l}(1-q^{m-k-i+\frac{1}{2}})} X_n^l.
$$
Comparing the coefficients of $X_n^l$,
we see that $(c_{k,m}(q))_{1\leq k,m\leq n}$ and \eqref{eq:invq} are indeed inverses.

\item Equation \eqref{eq:diff} may be written as
\begin{align}
(-1)^{n-l}X_l^m
-(-1)^{n-l+1}X_{l-1}^m =(-1)^{n-l}\sum_{k=0}^{\lfloor m/2 \rfloor}
g_{m,m-k}(q)\frac{(1-q^l)^{2m-2k}}{(1-q)^{2m-2k}}q^{-l(m-k)} \label{eq:diff2}
\end{align}
Summing Equation \eqref{eq:diff2} over $l$ from $1$ to $n$ and applying
Equation \eqref{eq:altersum}, we obtain
\begin{align}
X_n^m
=\sum_{k=0}^{\lfloor m/2 \rfloor}
g_{m,m-k}(q)T_{2m-2k,n}(q)q^{-n(m-k)}. \label{eq:gmntmn}
\end{align}
Substituting \eqref{eq:t2mnq} into \eqref{eq:gmntmn}, the right-hand side becomes
\begin{align}
\sum_{k=0}^{\lfloor m/2 \rfloor}\sum_{l=1}^{m-k}
g_{m,m-k}(q)(-1)^{m-k-l}
\frac{G_{m-k,m-k-l}(q)}{\prod_{i=0}^{m-k-l}(1+q^{m-k-i})}
X_n^l.
\end{align}
Comparing the coefficients of $X_n^l$,
we see that $(g_{k,m}(q))_{1\leq k,m \leq n}$ and \eqref{eq:invg} are inverse to each other.

\item Equation~\eqref{eq:sumd} may be written as
\begin{align}
 &\hskip -2mm  (-1)^{n-l}\frac{1-q^{l+\frac12}}{(1-q^\frac12)q^\frac{l}{2}}
X_l^{m-1}
     -(-1)^{n-l+1}\frac{1-q^{l-\frac12}}{(1-q^\frac12)q^\frac{l-1}{2}}
X_{l-1}^{m-1} \nonumber\\
     &=(-1)^{n-l}\sum_k d_{m,m-k}(q^\frac12)[l]^{2(m-k-\frac12)}q^{-l(m-k-\frac12)}. \label{eq:sumd2}
    \end{align}
Summing Equation \eqref{eq:sumd2} over $l$ from $1$ to $n$ and applying Equation
\eqref{eq:altersum}, we obtain
\begin{align}\label{eq:sumd3}
   \frac{1-q^{n+\frac12}}{(1-q^\frac12)q^\frac{n}{2}}X_n^{m-1}
    =\sum_k d_{m,m-k}(q^\frac12) T_{2m-2k-1,n}(q)q^{-n(m-k-\frac12)},\quad m\geq 2.
    \end{align}
Substituting \eqref{eq:t2m+1} into \eqref{eq:sumd3} yields
\begin{align}
    & \frac{1-q^{n+\frac12}}{(1-q^\frac12)q^\frac{n}{2}}\left(X_n^{m-1}
        -
         \sum_k \sum_{l=1}^{m-k}
          \frac{(-1)^{m-k-l}d_{m,m-k}(q^{\frac12})H_{m-k,m-k-l}(q^{\frac12})X_n^{l-1}}
                     {(1+q^{\frac12})^{m-k-l+1}\prod_{i=0}^{m-k-l}(1+q^{m-k-i-\frac12})}\right)
\nonumber\\
&=(-1)^n\sum_k \frac{(-1)^{m-k}d_{m,m-k}(q^{\frac12}) H_{m-k,m-k-1}(q^{\frac12})}
                     {(1+q^{\frac12})^{m-k}\prod_{i=0}^{m-k-1}(1+q^{m-k-i-\frac12})}.\label{eq:d-str}
\end{align}
We now show that the right-hand side of \eqref{eq:d-str} must vanish.
Suppose $0<q<1$. Denote the left-hand side of \eqref{eq:d-str} by $L_n$.
If there exists an $n\in\mathbb{N}$ such that $L_n=0$ we are done. Suppose $L_n\neq 0$ for all $n\geq 1$,
then  $L_n$ is  a rational function
in $t=q^{\frac n2}$ and can be written as
$$
L_n=t^sf(t) \qquad \textrm{with}\quad t=q^{\frac n2},
$$
where $s$ is an integer and $f(t)$ a rational function with $f(0)\neq 0$.
Since $f(q^{\frac n 2})\neq 0$, the right-hand side of \eqref{eq:d-str}
implies that
$$
f(q^{\frac n 2})f(q^{\frac{n+1}2})<0\qquad \forall\, n\geq 1.
$$
Taking the limit as $n\to \infty$ we get
$(f(0))^2\leq 0$, which is impossible.  Hence $L_n=0$ and \eqref{eq:d-str} reduces to
\begin{align}\label{eq:sumd4}
X_n^{m-1}
=\sum_k d_{m,m-k}(q^{\frac12})\sum_{l=1}^{m-k}
    \frac{(-1)^{m-k-l} H_{m-k,m-k-l}(q^{\frac12})X_n^{l-1}}
                     {(1+q^{\frac12})^{m-k-l+1}\prod_{i=0}^{m-k-l}(1+q^{m-k-i-\frac12})}.
    \end{align}
Comparing the coefficients of $X_n^{l-1}$ on
both sides of \eqref{eq:sumd4}, we see that $(d_{k,m}(q))_{1\leq k,m \leq n}$
and \eqref{eq:invh} are indeed inverses.
\end{itemize}
\end{proof}

The following easily verified result has been given by Gessel and
 Viennot~\cite{GesselViennot1989}.
\begin{lem}\label{lem:detinv}
  Let $(A_{ij})_{0\leq i,j \leq m}$ be an invertible lower triangular matrix,
  and let $(B_{ij})=(A_{ij})^{-1}$. Then for $0\leq k\leq n\leq m$, we have
  \begin{equation*}
     B_{n,k}=\frac{(-1)^{n-k}}{A_{k,k}A_{k+1,k+1}\cdots A_{n,n}}
             \left|A_{k+i+1,k+j}\right|_{0\leq i,j \leq n-k-1}.
  \end{equation*}
\end{lem}
Using the above lemma we derive immediately  from
Theorem~\ref{thm:inverse} the following determinant formulas:
\begin{align}
P_{m,k}(q)&=\det_{0\leq i,j \leq k-1}(h_{m-k-i+2j-1}(\{1,
q\}^{i-j+2})),\label{eq:detp}\\
Q_{m,k}(q)&=\det_{0\leq i,j \leq k-1}(c_{m-k+i+1,m-k+j}(q)),\label{eq:detq}\\
G_{m,k}(q) &=\det_{0\leq i,j\leq k-1}(g_{m-k+i+1,m-k+j}(q)),\label{eq:detg}\\
H_{m,k}(q) &=\det_{0\leq i,j\leq
k-1}(d_{m-k+i+1,m-k+j}(q)).\label{eq:deth}
\end{align}

\section{Combinatorial interpretations}
A \emph{lattice path} or \emph{path} $s_0\to s_n$ is a sequence of
points $(s_0,s_1,\ldots,s_n)$ in the plane $\mathbb Z^2$ such that
$s_i-s_{i-1}=(1,0),\,(0,1)$ for all $i=1,\ldots, n$. Let us assign
a weight to each step $(s_i,s_{i+1})$ of $s_0\to s_n$. We define
the weight $N(s_0\to s_n)$ of the path $s_0\to s_n$ to be the
product of the weights of its steps. Let $s_0=(a,b)$ and
$s_n=(c,d)$, if we weight each vertical step with $x$-coordinate
$i$ by $x_i$  and all horizontal steps by 1 then
\begin{equation}\label{eq:weight}
N(s_0\to s_n)=h_{d-b}(x_a, x_{a+1},\ldots, x_c).
\end{equation}
Now consider  two sequences of lattice points  ${\bf u}:=(u_1,\,
u_2,\, \ldots, u_n)$ and ${\bf v}:=(v_1,\, v_2,\, \ldots, v_n)$
such that for $i<j$ and $k<l$ any lattice path between $u_i$ and
$v_l$ has a common point with any lattice path between $u_j$ and
$v_k$. Set
$$
N({\bf u}, {\bf v}):=\sum N(u_1\to  v_1)\cdots N(u_n\to v_n),
$$
where the sum is over all families of non-intersecting paths
$(u_1\to v_1,\ldots,u_n\to v_n)$.

The following remarkable result can be found in Gessel and
Viennot~\cite{GesselViennot1989}. For historical remarks see also
Krattenthaler~\cite{Krattenthaler2005}.
\begin{thm}\label{thm:lgv}
[Lindstr\"om-Gessel-Viennot] We have
  $$
  N({\bf u},{\bf v})={\det_{1\leq i,j\leq n}(N(u_j\to v_i))}.
  $$
\end{thm}

We are now ready to exhibit the combinatorial interpretation of
the $q$-Faulhaber numbers.
\begin{thm}\label{thm:pq}
  Let ${\bf u}=(u_0, \ldots, u_{k-1})$ and ${\bf v}=(v_0, \ldots, v_{k-1})$, where
  $u_i=(2i, -2i)$ and  $v_i=(2i+3, m-k-i-1)$ for $0\leq i\leq k-1$.
  \begin{enumerate}
  \item The polynomial
    $P_{m,k}(q)$
    is the sum of the weights of $k$-non-intersecting paths
    from ${\bf u}$ to ${\bf v}$,
    where a vertical step with an even $x$-coordinate has weight $q$, and all
    the other steps have weight $1$.
  \item The polynomial
    $Q_{m,k}(q)$
    is the sum of the weights of $k$-non-intersecting paths
    from ${\bf u}$ to ${\bf v}$, where the weight of the individual steps is the same as
    before with the exception that $q$ is replaced with $q^2$
    and the vertical step starting from any $u_j$ has  weight $q^2+q$ instead of $q^2$.
  \end{enumerate}
\end{thm}
\begin{proof}
  For (i), by means of \eqref{eq:weight} we have
  $$
  N(u_j\to v_i)=h_{m-k-i+2j-1}(\{1, q\}^{i-j+2}).
  $$
The result then follows from \eqref{eq:detp} and Theorem~\ref{thm:lgv}.

  For (ii), assume that $u_{j}'=(2j+1,-2j)$ and  $u_{j}''=(2j,1-2j)$.
  The first step of a lattice path from $u_j$ to $v_i$ is
  either $u_j\to u_j'$ or $u_j\to u_j''$. As  $N(u_j\to u_j')=1$, $N(u_j\to u_j'')=q^2+q$ and
  $h_n(x_1,\ldots,x_{r-1})+x_r h_{n-1}(x_1,\ldots,x_{r})=h_n(x_1,\ldots,x_{r})$, we have
  \begin{align*}
    N(u_j\to v_i)
    &=N(u_j\to u_j')N(u_j'\to v_i)+N(u_j\to u_j'')N(u_j''\to v_i) \\
    &=N(u_j'\to v_i)+(q^2+q)N(u_j''\to v_i) \\
    &=h_{m-k-i+2j-1}(\{1\}^{i-j+2},\{q^2\}^{i-j+1})\\
    &\qquad \qquad +(q^2+q)h_{m-k-i+2j-2}(\{1,
    q^2\}^{i-j+2})
    \\
    &=h_{m-k-i+2j-1}(\{1, q^2\}^{i-j+2})+qh_{m-k-i+2j-2}(\{1, q^2\}^{i-j+2}).
  \end{align*}
The result then follows from \eqref{eq:detq} and Theorem~\ref{thm:lgv}.

\end{proof}

\begin{cor}
The polynomials $P_{m,k}(q)$ and $Q_{m,k}(q)$ have symmetric coefficients.
\end{cor}
\begin{proof}
  A combinatorial way to see the symmetry of the coefficients of $P_{m,k}(q)$
  is as follows: Modifying the weights such that vertical steps with an odd
  $x$-coordinate have weight $q$ and all the others weight $1$ does not change
  the entries of the determinant.

  However, consider any given family of paths with weight $q^w$, when vertical
  steps with even $x$-coordinate have weight $q$. After the modification of the
  weights it will have weight $q^{{\rm max}-w}$, where ${\rm max}$ is the total number of
  vertical steps in such a family of paths, which implies the claim.

  For the polynomials $Q_{m,k}$, we use the following alternative weight:
  vertical steps with odd $x$-coordinate have weight $q^2$, vertical steps with
  starting point $u_i$ have weight $1+q$ and all others have weight $1$.
\end{proof}

When $k=m-1$, there is only one lattice path from $u_0=(0,0)$ to $v_0=(3,0)$,
which has weight $1$. This establishes the following result:
\begin{cor}
For $m\geq 2$, we have $P_{m,m-1}(q)=P_{m,m-2}(q)$ and  $Q_{m,m-1}(q)=Q_{m,m-2}(q)$.
\end{cor}

For the combinatorial interpretation of the $q$-Sali\'e numbers, we need an
auxiliary lemma:

\begin{lem}\label{lem:detsum}
  Let $(A_{ij})_{1\leq i,j\leq n}$ and $(B_{ij})_{1\leq i,j\leq n}$ be two
  matrices. Then
  $$
  \det_{1\leq i,j\leq n}(A_{ij}+B_{ij})
  =\sum_{I\subseteq\{1,\ldots,n\}}\det_{1\leq i,j\leq n}(D_{ij}^{(I)}),
  $$
  where
  \begin{equation*}
 D_{ij}^{(I)}
    =\begin{cases}
       A_{ij}, &\text{if $j\in I$,} \\[5pt]
       B_{ij}, &\text{otherwise.}
     \end{cases}
  \end{equation*}
\end{lem}

\begin{thm}\label{thm:gh}
  Let ${\bf u}=(u_0, \ldots, u_{k-1})$ and ${\bf v}=(v_0, \ldots, v_{k-1})$, where
 $u_i=(2i, -2i)$ and $v_i=(2i+2, m-k-1-i)$ for $0\leq i\leq k-1$.
  \begin{itemize}
  \item[(i)] The polynomial
    $G_{m,k}(q)$
    is the sum of the weights of $k$-non-intersecting lattice paths
    ${\bf L}$ from ${\bf u}$ to ${\bf v}$
    with the weight of ${\bf L}$ being
    \begin{equation*}
      \sum_{I\subseteq\{0,1,\dots,k-1\}} w_I({\bf L}),
    \end{equation*}
    where $w_I$ is defined as follows: for each $i\in I$, vertical
    steps with $x$-coordinate $2i-1$ have weight $q$, and for any integer
    $i\notin I$, vertical steps with $x$-coordinate $2i$ have weight $q$.
    All other steps have weight $1$.

  \item[(ii)] The polynomial $H_{m,k}(q)$
    is the sum of the weights of $k$-non-intersecting lattice paths ${\bf L}$
    from ${\bf u}$ to ${\bf v}$,
    with the weight of ${\bf L}$ being
   \begin{equation*}
      \sum_{I\subseteq\{0,1,\dots,k-1\}} \overline{w}_I({\bf L}),
    \end{equation*}
    where $\overline{w}_I$ is the same as $w_I$ -- replacing $q$ with $q^2$ --
    with the exception of vertical steps starting from one of the points $u_i$,
    which have an additional weight of $q$. More precisely, if the weight of
    such a step would be $1$, it has weight $1+q$, it its weight would be
    $q^2$, it has weight $q^2+q$.
  \end{itemize}
\end{thm}

\begin{proof}
  (i) We apply Lemma~\ref{lem:detsum} to $\det_{0\leq i,j\leq
    k-1}(g_{m-k+i+1,m-k+j}(q))$, where
  \begin{equation*}
    \begin{split}
       g_{m-k+i+1,m-k+j}(q)&=h_{m-k-i+2j-1}(\{1\}^{i-j+2},\{q\}^{i-j+1})\\
       &\quad+h_{m-k-i+2j-1}(\{1\}^{i-j+1},\{q\}^{i-j+2}).
    \end{split}
  \end{equation*}
  Suppose that $j\in I$ and $0\leq i\leq k-1$. Then we have to show that
  $h_{m-k-i+2j-1}(\{1\}^{i-j+2},\{q\}^{i-j+1})$ is the weighted sum of lattice
  paths from $u_j$ to $v_i$, where the vertical steps have the weight given in
  the claim. To this end, note that
  $h_{m-k-i+2j-1}(\{1\}^{i-j+2},\{q\}^{i-j+1})$ counts lattice paths from $u_j$
  to $v_i$, when steps on $i-j+1$ given vertical lines have weight $q$, those
  steps on the remaining $i-j+2$ vertical lines have weight $1$.

  By the construction in the claim, steps on exactly one of the vertical lines
  with $x$-coordinates $2r-1$ and $2r$ have weight $q$. Since $j\in I$, steps
  on the vertical line with $x$-coordinate $2j$, i.e., with the
  $x$-coordinate of $u_j$, have weight $1$.

  Similarly, if $j\not\in I$ we can verify that there are exactly $i-j+2$ vertical
  lines between $u_j$ and $v_i$ with steps thereon having weight $q$.

(ii) In the same way, we can show that for $j\in I$ and $0\leq i\leq k-1$.
$$h_{m-k-i+2j-1}(\{1\}^{i-j+2},\{q^2\}^{i-j+1})+qh_{m-k-i+2j-2}(\{1\}^{i-j+2},\{q^2\}^{i-j+1})$$
 is the sum of weights of lattice paths from $u_j$ to $v_i$, where the vertical steps have the weight given in
 the claim. Meanwhile, for $j\notin I$ and $0\leq i\leq k-1$,
$$h_{m-k-i+2j-1}(\{1\}^{i-j+1},\{q^2\}^{i-j+2})+qh_{m-k-i+2j-2}(\{1\}^{i-j+1},\{q^2\}^{i-j+2})$$
 is the sum of weights of lattice paths from $u_j$ to $v_i$.
\end{proof}

As an illustration of the underlying configurations in Theorem
\ref{thm:gh}, we give an example in Figure~\ref{fig:wi} for $m=7$ and $k=4$.
\setlength{\unitlength}{1.4mm}
\begin{figure}[h!]
\caption{Example for $w_I$ in Theorem \ref{thm:gh}.\label{fig:wi}}
\begin{picture}(55,58)
\put(0,40){\vector(1,0){55}}
\put(5,10){\vector(0,1){45}}
\put(5,40){\circle*{1}} \put(3,38){\footnotesize $0$}
\put(15,30){\circle*{1}}\put(15,30){\circle{2}} \put(12,27){\footnotesize $(2,-2)$}
\put(25,20){\circle*{1}}\put(25,20){\circle{2}} \put(22,17){\footnotesize $(4,-4)$}
\put(35,10){\circle*{1}}\put(32,7){\footnotesize $(6,-6)$}
\put(15,50){\circle*{1}}\put(13,52){\footnotesize $(2,2)$}
\put(25,45){\circle*{1}}\put(23,47){\footnotesize $(4,1)$}
\put(35,40){\circle*{1}}\put(33,42){\footnotesize $(6,0)$}
\put(45,35){\circle*{1}}\put(43,37){\footnotesize $(8,-1)$}
\put(5,40){\textcolor{magenta}{\thicklines\vector(0,1){5}}}\put(5,45){\circle*{0.7}}
\put(5,45){\thicklines\vector(1,0){5}}\put(10,45){\circle*{0.7}}
\put(10,45){\textcolor{magenta}{\thicklines\vector(0,1){5}}}\put(10,50){\circle*{0.7}}
\put(10,50){\thicklines\vector(1,0){5}}
\put(15,30){\thicklines\vector(0,1){10}}\put(15,35){\circle*{0.7}}\put(15,40){\circle*{0.7}}
\put(15,40){\thicklines\vector(1,0){5}}\put(20,40){\circle*{0.7}}
\put(20,40){\textcolor{magenta}{\thicklines\vector(0,1){5}}}\put(20,45){\circle*{0.7}}
\put(20,45){\thicklines\vector(1,0){5}}
\put(25,20){\thicklines\vector(0,1){5}}\put(25,25){\circle*{0.7}}
\put(25,25){\thicklines\vector(1,0){5}}\put(30,25){\circle*{0.7}}
\put(30,25){\thicklines\vector(0,1){5}}\put(30,30){\circle*{0.7}}
\put(30,30){\thicklines\vector(1,0){5}}\put(35,30){\circle*{0.7}}
\put(35,30){\textcolor{magenta}{\thicklines\vector(0,1){10}}}\put(35,35){\circle*{0.7}}
\put(35,10){\thicklines\vector(1,0){5}} \put(40,10){\circle*{0.7}}
\put(40,10){\thicklines\vector(0,1){10}}\put(40,15){\circle*{0.7}}\put(40,20){\circle*{0.7}}
\put(40,20){\thicklines\vector(1,0){5}} \put(45,20){\circle*{0.7}}
\put(45,20){\textcolor{magenta}{\thicklines\vector(0,1){15}}}
\put(45,25){\circle*{0.7}}\put(45,30){\circle*{0.7}}
\put(0,2){$I=\{1,2\}$, $w_I({\bf L})=q^8$ and $\overline{w}_I({\bf L})=q^{14}(q+q^2)(q+1)^2$}
\end{picture}
\end{figure}

\begin{cor}
The polynomials $G_{m,k}(q)$ and $H_{m,k}(q)$ have symmetric coefficients.
\end{cor}
\begin{proof}
  A combinatorial way to see the symmetry of the coefficients of $G_{m,k}(q)$
  is as follows: Modifying $w_I$ such that for each $i\in I$, vertical steps
  with $x$-coordinate $2i$ have weight $q$, and for any integer $i\notin I$,
  vertical steps with $x$-coordinate $2i-1$ have weight $1$ does not change the
  entries of the determinant.

  However, consider any given family of paths with weight $q^w$ with weight by
  Theorem~\ref{thm:gh}(i). After the modification of the weights it will have
  weight $q^{{\rm max}-w}$, where ${\rm max}$ is the total number of vertical
  steps in such a family of paths, which implies the claim.

  We omit the proof of the symmetry of the coefficients of $H_{m,k}(q)$.
\end{proof}

\begin{cor} Let ${\bf u}=(u_0, \ldots, u_{k-1})$ and ${\bf v}=(v_0, \ldots, v_{k-1})$, where
 $u_i=(2i, -2i)$ and $v_i=(2i+2, m-k-1-i)$ for $0\leq i\leq k-1$.
\begin{itemize}
  \item[(i)] The polynomial $G_{m,k}(q)$
    is the sum of the weights of $k$-non-intersecting lattice paths
    ${\bf L}$ from ${\bf u}$ to ${\bf v}$ with the weight of ${\bf L}$ being
    \begin{equation*}
     q^{\sigma_{2k}({\bf L})}
    \prod_{i=0}^{k-1}\left(q^{\sigma_{2i-1}({\bf L})}+q^{\sigma_{2i}({\bf L})}\right),
    \end{equation*}
  where $\sigma_j$ denotes the number of vertical steps with $x$-coordinate $j$.

  \item[(ii)] The polynomial $H_{m,k}(q)$
    is the sum of the weights of $k$-non-intersecting lattice paths ${\bf L}$
    from ${\bf u}$ to ${\bf v}$ with the weight of ${\bf L}$ being
   \begin{equation*}
     (1+q)^{f({\bf L})}q^{2\sigma_{2k}({\bf L})}\prod_{i=0}^{k-1}
    \left(q^{2\sigma_{2i-1}({\bf L})}+q^{2\sigma_{2i}({\bf L})-f_i({\bf L})}\right),
    \end{equation*}
  where $\sigma_j$ is as in {\rm(i)} and $f$ {\rm(}resp. $f_i${\rm)} denotes the number
  of vertical steps starting from ${\bf u}$ {\rm(}resp. $u_i${\rm)}.
  \end{itemize}
\end{cor}

\begin{proof} (i) By the definition of $w_I$, for $0\leq i\leq k-1$, if $i\in I$, then
  vertical steps on the line with $x$-coordinates $2i-1$ have weight $q$
  and vertical steps on the line with $x$-coordinates $2i$ have weight $1$;
  and if $i\notin I$, the case is just contrary. Note that
  steps on the vertical line with $x$-coordinates $2k$ always have weight $q$
  and steps on the vertical line with $x$-coordinates $2k-1$ always have weight $1$.
  This implies that
    \begin{equation*}
      \sum_{I\subseteq\{0,1,\dots,k-1\}} w_I({\bf L})
      =q^{\sigma_{2k}({\bf L})}
    \prod_{i=0}^{k-1}\left(q^{\sigma_{2i-1}({\bf L})}+q^{\sigma_{2i}({\bf L})}\right).
    \end{equation*}

(ii) Notice that for $0\leq i\leq k-1$, we have $f_i({\bf L})=1$ if ${\bf L}$ contains
 a vertical step starting from the point $u_i$, and $f_i({\bf L})=0$ otherwise. Similarly, we have
\begin{align*}
   \sum_{I\subseteq\{0,1,\dots,k-1\}} \overline{w}_I({\bf L})
     & =q^{2\sigma_{2k}({\bf L})}
    \prod_{i=0}^{k-1}\left(q^{2\sigma_{2i-1}({\bf L})}(1+q)^{f_i({\bf L})}
      +q^{2\sigma_{2i}({\bf L})-2f_i({\bf L})}(q^2+q)^{f_i({\bf L})}\right),\\
    & =(1+q)^{f({\bf L})}q^{2\sigma_{2k}({\bf L})}\prod_{i=0}^{k-1}
    \left(q^{2\sigma_{2i-1}({\bf L})}+q^{2\sigma_{2i}({\bf L})-f_i({\bf L})}\right).
    \end{align*}
This completes the proof.
\end{proof}

\setlength{\unitlength}{1.6mm}
\begin{figure}[!h]
\caption{An illustration for
$G_{4,2}(q)=10q^3+24q^2+24q+10$.\label{fig:g}}
$$\begin{array}{cccc}
\begin{picture}(17,16)
\put(0,4){\vector(0,1){12}} \put(0,11){\vector(1,0){15}}
\put(0,11){\circle*{0.75}} \put(-1.5,10){\footnotesize $0$}
\put(6,14){\circle*{0.75}} \put(5.3,15){\footnotesize $(2,1)$}
\put(6,5){\circle*{0.75}} \put(4,3.2){\footnotesize $(2,-2)$}
\put(12,11){\circle*{0.75}} \put(12.6,9.5){\footnotesize $4$}
\put(0,11){\thicklines\vector(0,1){3}}
\put(0,14){\thicklines\vector(1,0){6}}
\put(6,5){\thicklines\vector(0,1){6}}
\put(6,11){\thicklines\vector(1,0){6}} \put(3,0.5){\footnotesize
$q^3+q^2+q+1$}
\end{picture}
&
\begin{picture}(17,16)
\put(0,4){\vector(0,1){12}} \put(0,11){\vector(1,0){15}}
\put(0,11){\circle*{0.75}} \put(-1.5,10){\footnotesize $0$}
\put(6,14){\circle*{0.75}} \put(5.3,15){\footnotesize $(2,1)$}
\put(6,5){\circle*{0.75}} \put(4,3.2){\footnotesize $(2,-2)$}
\put(12,11){\circle*{0.75}} \put(12.6,9.5){\footnotesize $4$}
\put(0,11){\thicklines\vector(1,0){3}}
\put(3,11){\thicklines\vector(0,1){3}}
\put(3,14){\thicklines\vector(1,0){3}}
\put(6,5){\thicklines\vector(0,1){6}}
\put(6,11){\thicklines\vector(1,0){6}} \put(3,0.5){\footnotesize
$2q^2+2q$}
\end{picture}
&
\begin{picture}(17,16)
\put(0,4){\vector(0,1){12}} \put(0,11){\vector(1,0){15}}
\put(0,11){\circle*{0.75}} \put(-1.5,10){\footnotesize $0$}
\put(6,14){\circle*{0.75}} \put(5.3,15){\footnotesize $(2,1)$}
\put(6,5){\circle*{0.75}} \put(4,3.2){\footnotesize $(2,-2)$}
\put(12,11){\circle*{0.75}} \put(12.6,9.5){\footnotesize $4$}
\put(0,11){\thicklines\vector(0,1){3}}
\put(0,14){\thicklines\vector(1,0){6}}
\put(6,5){\thicklines\vector(0,1){3}}
\put(6,8){\thicklines\vector(1,0){3}}
\put(9,8){\thicklines\vector(0,1){3}}
\put(9,11){\thicklines\vector(1,0){3}} \put(3,0.5){\footnotesize
$q^2+2q+1$}
\end{picture}
&
\begin{picture}(17,16)
\put(0,4){\vector(0,1){12}} \put(0,11){\vector(1,0){15}}
\put(0,11){\circle*{0.75}} \put(-1.5,10){\footnotesize $0$}
\put(6,14){\circle*{0.75}} \put(5.3,15){\footnotesize $(2,1)$}
\put(6,5){\circle*{0.75}} \put(4,3.2){\footnotesize $(2,-2)$}
\put(12,11){\circle*{0.75}} \put(12.6,9.5){\footnotesize $4$}
\put(0,11){\thicklines\vector(1,0){3}}
\put(3,11){\thicklines\vector(0,1){3}}
\put(3,14){\thicklines\vector(1,0){3}}
\put(6,5){\thicklines\vector(0,1){3}}
\put(6,8){\thicklines\vector(1,0){3}}
\put(9,8){\thicklines\vector(0,1){3}}
\put(9,11){\thicklines\vector(1,0){3}} \put(3,0.5){\footnotesize
$4q$}
\end{picture} \\[15pt]
\begin{picture}(17,16)
\put(0,4){\vector(0,1){12}} \put(0,11){\vector(1,0){15}}
\put(0,11){\circle*{0.75}} \put(-1.5,10){\footnotesize $0$}
\put(6,14){\circle*{0.75}} \put(5.3,15){\footnotesize $(2,1)$}
\put(6,5){\circle*{0.75}} \put(4,3.2){\footnotesize $(2,-2)$}
\put(12,11){\circle*{0.75}} \put(12.6,9.5){\footnotesize $4$}
\put(0,11){\thicklines\vector(1,0){6}}
\put(6,11){\thicklines\vector(0,1){3}}
\put(6,5){\thicklines\vector(0,1){3}}
\put(6,8){\thicklines\vector(1,0){3}}
\put(9,8){\thicklines\vector(0,1){3}}
\put(9,11){\thicklines\vector(1,0){3}} \put(3,0.5){\footnotesize
$2q^2+2$}
\end{picture}
&
\begin{picture}(17,16)
\put(0,4){\vector(0,1){12}} \put(0,11){\vector(1,0){15}}
\put(0,11){\circle*{0.75}} \put(-1.5,10){\footnotesize $0$}
\put(6,14){\circle*{0.75}} \put(5.3,15){\footnotesize $(2,1)$}
\put(6,5){\circle*{0.75}} \put(4,3.2){\footnotesize $(2,-2)$}
\put(12,11){\circle*{0.75}} \put(12.6,9.5){\footnotesize $4$}
\put(0,11){\thicklines\vector(0,1){3}}
\put(0,14){\thicklines\vector(1,0){6}}
\put(6,5){\thicklines\vector(0,1){3}}
\put(6,8){\thicklines\vector(1,0){6}}
\put(12,8){\thicklines\vector(0,1){3}} \put(3,0.5){\footnotesize
$q^3+2q^2+q$}
\end{picture}
&
\begin{picture}(17,16)
\put(0,4){\vector(0,1){12}} \put(0,11){\vector(1,0){15}}
\put(0,11){\circle*{0.75}} \put(-1.5,10){\footnotesize $0$}
\put(6,14){\circle*{0.75}} \put(5.3,15){\footnotesize $(2,1)$}
\put(6,5){\circle*{0.75}} \put(4,3.2){\footnotesize $(2,-2)$}
\put(12,11){\circle*{0.75}} \put(12.6,9.5){\footnotesize $4$}
\put(0,11){\thicklines\vector(1,0){3}}
\put(3,11){\thicklines\vector(0,1){3}}
\put(3,14){\thicklines\vector(1,0){3}}
\put(6,5){\thicklines\vector(0,1){3}}
\put(6,8){\thicklines\vector(1,0){6}}
\put(12,8){\thicklines\vector(0,1){3}} \put(3,0.5){\footnotesize
$4q^2$}
\end{picture}
&
\begin{picture}(17,16)
\put(0,4){\vector(0,1){12}} \put(0,11){\vector(1,0){15}}
\put(0,11){\circle*{0.75}} \put(-1.5,10){\footnotesize $0$}
\put(6,14){\circle*{0.75}} \put(5.3,15){\footnotesize $(2,1)$}
\put(6,5){\circle*{0.75}} \put(4,3.2){\footnotesize $(2,-2)$}
\put(12,11){\circle*{0.75}} \put(12.6,9.5){\footnotesize $4$}
\put(0,11){\thicklines\vector(1,0){6}}
\put(6,11){\thicklines\vector(0,1){3}}
\put(6,5){\thicklines\vector(0,1){3}}
\put(6,8){\thicklines\vector(1,0){6}}
\put(12,8){\thicklines\vector(0,1){3}} \put(3,0.5){\footnotesize
$2q^3+2q$}
\end{picture} \\[15pt]

\begin{picture}(17,16)
\put(0,4){\vector(0,1){12}} \put(0,11){\vector(1,0){15}}
\put(0,11){\circle*{0.75}} \put(-1.5,10){\footnotesize $0$}
\put(6,14){\circle*{0.75}} \put(5.3,15){\footnotesize $(2,1)$}
\put(6,5){\circle*{0.75}} \put(4,3.2){\footnotesize $(2,-2)$}
\put(12,11){\circle*{0.75}} \put(12.6,9.5){\footnotesize $4$}
\put(0,11){\thicklines\vector(0,1){3}}
\put(0,14){\thicklines\vector(1,0){6}}
\put(6,5){\thicklines\vector(1,0){3}}
\put(9,5){\thicklines\vector(0,1){6}}
\put(9,11){\thicklines\vector(1,0){3}} \put(3,0.5){\footnotesize
$2q+2$}
\end{picture}
&
\begin{picture}(17,16)
\put(0,4){\vector(0,1){12}} \put(0,11){\vector(1,0){15}}
\put(0,11){\circle*{0.75}} \put(-1.5,10){\footnotesize $0$}
\put(6,14){\circle*{0.75}} \put(5.3,15){\footnotesize $(2,1)$}
\put(6,5){\circle*{0.75}} \put(4,3.2){\footnotesize $(2,-2)$}
\put(12,11){\circle*{0.75}} \put(12.6,9.5){\footnotesize $4$}
\put(0,11){\thicklines\vector(1,0){3}}
\put(3,11){\thicklines\vector(0,1){3}}
\put(3,14){\thicklines\vector(1,0){3}}
\put(6,5){\thicklines\vector(1,0){3}}
\put(9,5){\thicklines\vector(0,1){6}}
\put(9,11){\thicklines\vector(1,0){3}} \put(3,0.5){\footnotesize
$2q+2$}
\end{picture}
&
\begin{picture}(17,16)
\put(0,4){\vector(0,1){12}} \put(0,11){\vector(1,0){15}}
\put(0,11){\circle*{0.75}} \put(-1.5,10){\footnotesize $0$}
\put(6,14){\circle*{0.75}} \put(5.3,15){\footnotesize $(2,1)$}
\put(6,5){\circle*{0.75}} \put(4,3.2){\footnotesize $(2,-2)$}
\put(12,11){\circle*{0.75}} \put(12.6,9.5){\footnotesize $4$}
\put(0,11){\thicklines\vector(1,0){6}}
\put(6,11){\thicklines\vector(0,1){3}}
\put(6,5){\thicklines\vector(1,0){3}}
\put(9,5){\thicklines\vector(0,1){6}}
\put(9,11){\thicklines\vector(1,0){3}} \put(3,0.5){\footnotesize
$2q+2$}
\end{picture}
&
\begin{picture}(17,16)
\put(0,4){\vector(0,1){12}} \put(0,11){\vector(1,0){15}}
\put(0,11){\circle*{0.75}} \put(-1.5,10){\footnotesize $0$}
\put(6,14){\circle*{0.75}} \put(5.3,15){\footnotesize $(2,1)$}
\put(6,5){\circle*{0.75}} \put(4,3.2){\footnotesize $(2,-2)$}
\put(12,11){\circle*{0.75}} \put(12.6,9.5){\footnotesize $4$}
\put(0,11){\thicklines\vector(0,1){3}}
\put(0,14){\thicklines\vector(1,0){6}}
\put(6,5){\thicklines\vector(1,0){3}}
\put(9,5){\thicklines\vector(0,1){3}}
\put(9,8){\thicklines\vector(1,0){3}}
\put(12,8){\thicklines\vector(0,1){3}} \put(3,0.5){\footnotesize
$2q^2+2q$}
\end{picture}
\\[15pt]
\begin{picture}(17,16)
\put(0,4){\vector(0,1){12}} \put(0,11){\vector(1,0){15}}
\put(0,11){\circle*{0.75}} \put(-1.5,10){\footnotesize $0$}
\put(6,14){\circle*{0.75}} \put(5.3,15){\footnotesize $(2,1)$}
\put(6,5){\circle*{0.75}} \put(4,3.2){\footnotesize $(2,-2)$}
\put(12,11){\circle*{0.75}} \put(12.6,9.5){\footnotesize $4$}
\put(0,11){\thicklines\vector(1,0){3}}
\put(3,11){\thicklines\vector(0,1){3}}
\put(3,14){\thicklines\vector(1,0){3}}
\put(6,5){\thicklines\vector(1,0){3}}
\put(9,5){\thicklines\vector(0,1){3}}
\put(9,8){\thicklines\vector(1,0){3}}
\put(12,8){\thicklines\vector(0,1){3}} \put(3,0.5){\footnotesize
$2q^2+2q$}
\end{picture}
&
\begin{picture}(17,16)
\put(0,4){\vector(0,1){12}} \put(0,11){\vector(1,0){15}}
\put(0,11){\circle*{0.75}} \put(-1.5,10){\footnotesize $0$}
\put(6,14){\circle*{0.75}} \put(5.3,15){\footnotesize $(2,1)$}
\put(6,5){\circle*{0.75}} \put(4,3.2){\footnotesize $(2,-2)$}
\put(12,11){\circle*{0.75}} \put(12.6,9.5){\footnotesize $4$}
\put(0,11){\thicklines\vector(1,0){6}}
\put(6,11){\thicklines\vector(0,1){3}}
\put(6,5){\thicklines\vector(1,0){3}}
\put(9,5){\thicklines\vector(0,1){3}}
\put(9,8){\thicklines\vector(1,0){3}}
\put(12,8){\thicklines\vector(0,1){3}} \put(3,0.5){\footnotesize
$2q^2+2q$}
\end{picture}
&
\begin{picture}(17,16)
\put(0,4){\vector(0,1){12}} \put(0,11){\vector(1,0){15}}
\put(0,11){\circle*{0.75}} \put(-1.5,10){\footnotesize $0$}
\put(6,14){\circle*{0.75}} \put(5.3,15){\footnotesize $(2,1)$}
\put(6,5){\circle*{0.75}} \put(4,3.2){\footnotesize $(2,-2)$}
\put(12,11){\circle*{0.75}} \put(12.6,9.5){\footnotesize $4$}
\put(0,11){\thicklines\vector(0,1){3}}
\put(0,14){\thicklines\vector(1,0){6}}
\put(6,5){\thicklines\vector(1,0){6}}
\put(12,5){\thicklines\vector(0,1){6}} \put(3,0.5){\footnotesize
$2q^3+2q^2$}
\end{picture}
&
\begin{picture}(17,16)
\put(0,4){\vector(0,1){12}} \put(0,11){\vector(1,0){15}}
\put(0,11){\circle*{0.75}} \put(-1.5,10){\footnotesize $0$}
\put(6,14){\circle*{0.75}} \put(5.3,15){\footnotesize $(2,1)$}
\put(6,5){\circle*{0.75}} \put(4,3.2){\footnotesize $(2,-2)$}
\put(12,11){\circle*{0.75}} \put(12.6,9.5){\footnotesize $4$}
\put(0,11){\thicklines\vector(1,0){3}}
\put(3,11){\thicklines\vector(0,1){3}}
\put(3,14){\thicklines\vector(1,0){3}}
\put(6,5){\thicklines\vector(1,0){6}}
\put(12,5){\thicklines\vector(0,1){6}} \put(3,0.5){\footnotesize
$2q^3+2q^2$}
\end{picture}
\\[15pt]
\begin{picture}(17,16)
\put(0,4){\vector(0,1){12}} \put(0,11){\vector(1,0){15}}
\put(0,11){\circle*{0.75}} \put(-1.5,10){\footnotesize $0$}
\put(6,14){\circle*{0.75}} \put(5.3,15){\footnotesize $(2,1)$}
\put(6,5){\circle*{0.75}} \put(4,3.2){\footnotesize $(2,-2)$}
\put(12,11){\circle*{0.75}} \put(12.6,9.5){\footnotesize $4$}
\put(0,11){\thicklines\vector(1,0){6}}
\put(6,11){\thicklines\vector(0,1){3}}
\put(6,5){\thicklines\vector(1,0){6}}
\put(12,5){\thicklines\vector(0,1){6}} \put(3,0.5){\footnotesize
$2q^3+2q^2$}
\end{picture}&  & &
\end{array}
$$
\end{figure}
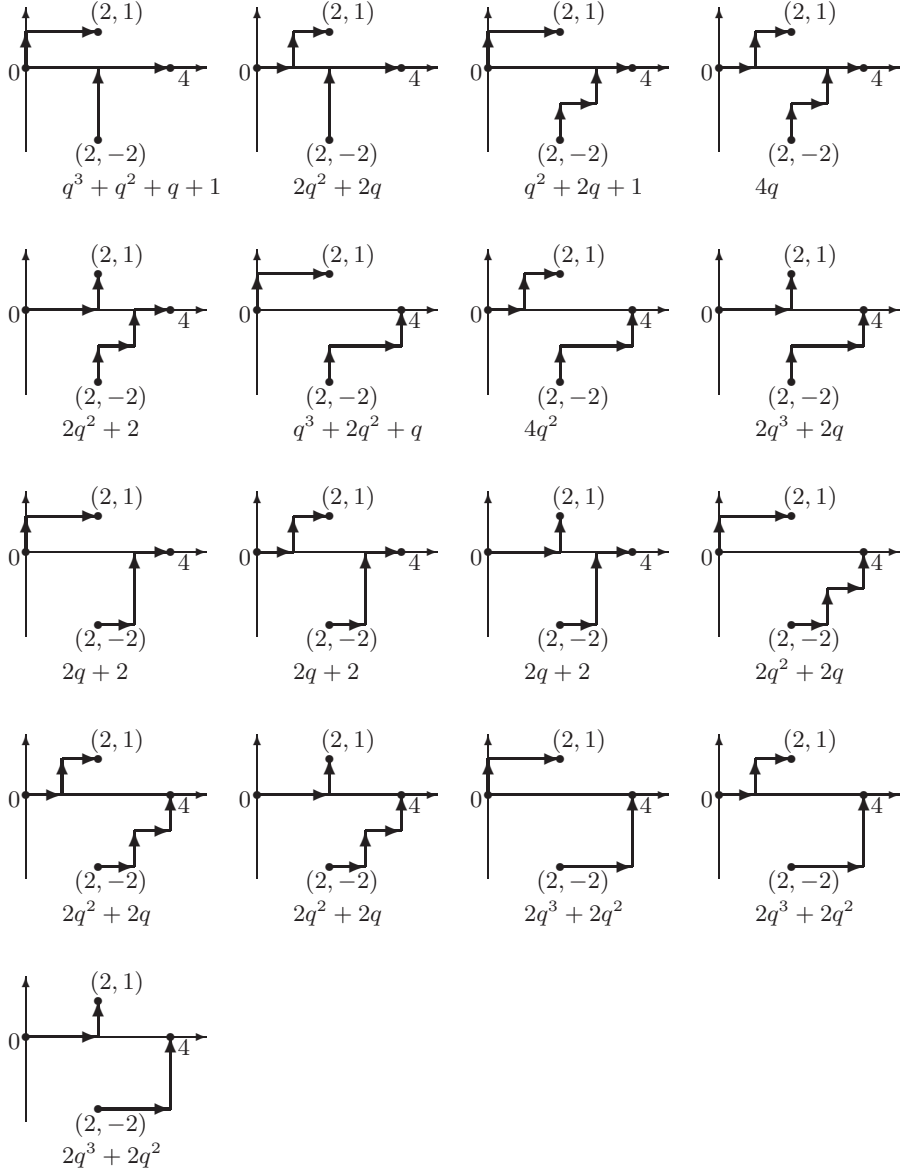

The computation of $G_{4,2}(q)$ is illustrated in Figure~\ref{fig:g},
while the value of $H_{4,2}(q)$ as given in Table~\ref{t:h} is computed in
Table~\ref{t:comh}.

\begin{table}[!h]
\caption{Values of $\sum_{I\subseteq \{0,1\}}\overline{w}_I({\bf L})$ corresponding
to Figure~\ref{fig:g}.\label{t:comh}}
{\footnotesize
\begin{center}
\begin{tabular}{|l|l|l|l|}
\hline
$(1+q)^3(1+q^3)$ & $2q^2(1+q)^2$ & $(1+q)^4$     & $2q(1+q)^2$ \\\hline
$2(1+q)(1+q^3)$  & $q^2(1+q)^4$  & $2q^3(1+q)^2$ & $2q^2(1+q)(1+q^3)$ \\\hline
$2(1+q)^2$       & $2(1+q^2)$    & $2(1+q^2)$    & $2q^2(1+q)^2$  \\\hline
$2q^2(1+q^2)$    & $2q^2(1+q^2)$ & $2q^4(1+q)^2$ & $2q^4(1+q^2)$  \\\hline
$2q^4(1+q^2)$    &               &               &   \\\hline
\end{tabular}
\end{center}
}
\end{table}

\begin{rmk}
Since
$$
\det_{1\leq i,j\leq n}(A_{ij}+B_{ij})
  =\sum_{I\subseteq\{1,\ldots,n\}}\det_{1\leq i,j\leq n}(C_{ij}^{(I)}),
$$
where
$$
 C_{ij}^{(I)}
    =\begin{cases}
       A_{ij}, &\text{if $i\in I$,} \\[5pt]
       B_{ij}, &\text{otherwise,}
     \end{cases}
$$
we may also define $w_I$ in Theorem~\ref{thm:gh}(i) as follows:
for each $i\in I$, vertical
    steps with $x$-coordinate $2i+3$ have weight $q$, and for any integer
    $i\notin I$, vertical steps with $x$-coordinate $2i+2$ have weight $q$.
    All other steps have weight $1$. In this case, for each
    $i\in I$ and $0\leq j\leq k-1$, we can show that
    $h_{m-k-i+2j-1}(\{1\}^{i-j+2},\{q\}^{i-j+1})$
    is the weighted sum of lattice paths from $u_j$ to $v_i$. Moreover,
     \begin{equation*}
      \sum_{I\subseteq\{0,1,\dots,k-1\}} w_I({\bf L})
      =q^{\sigma_{0}({\bf L})}
    \prod_{i=1}^{k}\left(q^{\sigma_{2i}({\bf L})}+q^{\sigma_{2i+1}({\bf L})}\right).
    \end{equation*}

 Similarly, we may define $\overline{w}_I$ in Theorem~\ref{thm:gh}(ii) as follows:
for each $i\in I$, a vertical step toward the point $v_i$ has weight $q+1$, vertical steps
with $x$-coordinate $2i+3$ have weight $q^2$. For any integer $i\notin I$,
a vertical step toward the point $v_i$ has weight $q^2+q$, and
vertical steps with $x$-coordinate $2i+2$ not toward $v_i$ have weight
$q^2$. All other steps have weight $1$. In this case, we have
\begin{equation*}
   \sum_{I\subseteq\{0,1,\dots,k-1\}} \overline{w}_I({\bf L})
   =(1+q)^{\overline{f}({\bf L})}q^{2\sigma_{0}({\bf L})}\prod_{i=1}^{k}
    \left(q^{2\sigma_{2i}({\bf L})-\overline{f}_i({\bf L})}
    +q^{2\sigma_{2i+1}({\bf L})}\right),
    \end{equation*}
 where $\overline{f}$ {\rm(}resp. $\overline{f}_i${\rm)} denotes the number
  of vertical steps ending in ${\bf v}$ {\rm(}resp. $v_i${\rm)}.

\end{rmk}

When $k=m-1$, there is only one lattice path from $u_0=(0,0)$ to $v_0=(2,0)$,
which has weight $1$. This establishes the following result:
\begin{cor}
  $G_{m,m-1}(q)=2G_{m,m-2}(q)$ and $H_{m,m-1}(q)=2H_{m,m-2}(q)$.
\end{cor}

\section{Open problems}

We would like to point out three directions of possible further research: It
appears that the polynomials $P_{m,k}$ and $G_{m,k}$ are log-concave, however,
we did not pursue this question further. Note that the polynomials $Q_{m,k}$
and $H_{m,k}$ are not even unimodal.

Victor Guo and Jiang Zeng gave in \cite{GuoZeng2005} even finer $q$-analogues
of the polynomials considered here, replacing \eqref{eq:qpower} and
\eqref{eq:altersum} by
\begin{align*}
  S_{m,n,r}(q) &=\sum_{k=1}^{n}\frac{[2rk]}{[2r]}
                             [k]^{m-1}q^{\frac{m+2r-1}{2}(n-k)},\\
  T_{m,n,r}(q) &=\sum_{k=1}^{n}(-1)^{n-k}
                             \frac{[(2r-1)k]}{[2r-1]}
                             [k]^{m-1} q^{\frac{m}{2}(n-k)},
\end{align*}
where $r\geq 1$.

Although the coefficients of the corresponding polynomials $P_{m,k,r},
Q_{m,k,r}, G_{m,k,r}$ and $H_{m,k,r}$ are not positive anymore, one might hope
for a refinement of Theorem~\ref{thm:inverse}.

Finally, we should point out that Ira Gessel and Xavier Viennot
\cite{GesselViennot1989} also presented nice generating functions for their
coefficients $f_{m,k}$ and $s_{m,k}$, namely
\begin{align*}
  \sum_{m,k}s_{m,k}t^k\frac{x^{2n}}{(2n)!}
  &=\frac{\cosh\sqrt{1+4t}\frac{x}{2}}{\cosh\frac{x}{2}},\\
  \sum_{m,k}f_{m,k}t^k\frac{x^{2n+1}}{(2n+1)!}
  &=\frac{\cosh\sqrt{1+4t}\frac{x}{2}-\cosh\frac{x}{2}}{t\sinh\frac{x}{2}}.
\end{align*}
It would be interesting to find the corresponding refinements.

\section{Epilogue}

One may wonder how these results were discovered. The truth is, that at first
\lq\lq only\rq\rq\ formula~\eqref{eq:p} was known. Using this formula,
Table~\ref{t:p} was computed. Then, in analogy to \cite{GesselViennot1989}, the
matrix
$$
\left((-1)^{k-m}\frac{[m]!}{[k+1]!}P_{k,k-m}(q)
\right)_{0\leq k,m\leq n}
$$
was inverted and, since we were looking for a lattice path
interpretation, the entry in row $i$ and column $j$ of the inverse
matrix had to be the weighted number of lattice paths from $u_j$
to $v_i$. This given, it was easy to find the correct weights.
Finally, we read the proof given in \cite{GesselViennot1989}
backwards, its first line corresponding to our
Lemma~\ref{lem:binomial}.

\vskip 5mm \noindent{\bf Acknowledgment.} The third  author was
supported by EC's IHRP Programme, within Research Training Network
``Algebraic Combinatorics in Europe," grant HPRN-CT-2001-00272.

\bibliographystyle{amsplain}

\end{document}